\documentclass[12pt]{amsart}
\usepackage{amsmath,amsthm}
\usepackage{amssymb}
\usepackage{euscript}			% nice script capitals (only capitals)
\usepackage{times}				% wonderfully improves on-screen fonts

\usepackage[matrix,arrow,curve,cmtip]{xy}

\xymatrixcolsep{1.9pc}                     % adjust size of diagrams
\xymatrixrowsep{1.9pc}
\newdir{ >}{{}*!/-5pt/\dir{>}}             % make better tailed arrows
\newdir{ |}{{}*!/-5pt/\dir{|}}             % make better mapsto arrows

\newtheorem{thm}{Theorem}[section]
\newtheorem{prop}[thm]{Proposition}

\newtheorem{lemma}[thm]{Lemma}

\theoremstyle{definition}

\newtheorem{conj}[thm]{Conjecture}

\theoremstyle{remark}

\numberwithin{equation}{section}

\DeclareMathOperator{\Ex}{Ex}

\DeclareMathOperator{\Sd}{Sd}

\DeclareMathOperator{\Sing}{Sing}

\newcommand{\cat}{\mathsf{cat}}
\newcommand{\Cat}{\mathit{Cat}}

\newcommand{\cof}{\mathsf{cof}}

\newcommand{\fib}{\mathsf{fib}}

\newcommand{\fin}{\mathsf{fin}}

\newcommand{\SSet}{\mathit{SSet}}

\renewcommand{\th}{{\text{th}}}
\newcommand{\Top}{\mathit{Top}}

\newcommand{\topos}{\EuScript}
\newcommand{\cate}{\mathcal}

\newcommand{\C}{\cate C}

\newcommand{\F}{\topos F}

\newcommand{\N}{\mathbb{N}}

\newcommand{\tW}{\mathsf W}

\newcommand{\lc}{\langle}					% left corner
\newcommand{\rc}{\rangle}					% right corner
\renewcommand{\leq}{\leqslant}				% the better-looking alternate
\renewcommand{\geq}{\geqslant}				% the better-looking alternate
\renewcommand{\subsetneq}{\varsubsetneqq}   % the better-looking alternate
\renewcommand{\supsetneq}{\supsetneqq}      % the better-looking alternate
\newcommand{\ra}{\rightarrow}
\newcommand{\lra}{\longrightarrow}

\newcommand{\inc}{\hookrightarrow}

\newcommand{\leri}{\leftrightarrows}

\newcommand{\msm}[1]{{\scriptstyle #1}}

\newcommand{\bd}{\textbf}

\newcommand{\three}[3]{\overset{#1}{\underset{#3}{#2}}}

\newcommand{\llra}[1]{\stackrel{#1}{\lra}}	% labeled long right arrow
\begin{document}

\title{Fibrations of simplicial sets}
\author{Tibor Beke}
\address{Department of Mathematics\\University of Massachusetts Lowell\\
One University Avenue\\Lowell, MA 01854}
\email{tibor\_beke@uml.edu}
\date{\today}
\begin{abstract}
There are infinitely many variants of the notion of Kan fibration that,
together with suitable choices of cofibrations and the usual notion of
weak equivalence of simplicial sets, satisfy Quillen's axioms for a homotopy model category.  
The combinatorics underlying these fibrations is purely finitary and seems
interesting both for its own sake and for its interaction with homotopy
types.  To show that these notions of fibration are indeed distinct, one needs to understand how iterates of Kan's Ex functor act on graphs and on nerves of small categories.
\end{abstract}
\maketitle

\section{Introduction}
The definition of fibration that now bears his name was introduced by Daniel
Kan in 1957, and remains a cornerstone of simplicial algebraic topology.
A decade later, Quillen axiomatized homotopy theory via his notion of a
\emph{model category} that comes equipped with three distinguished classes
of morphisms: fibrations, weak equivalences and cofibrations.  The
category of simplicial sets, where Kan fibrations, topological (also
called `combinatorial') weak equivalences, and monomorphisms serve these
roles, remains the primordial example of a homotopy model category.  The
goal of this article is to prove the following

\bd{Theorem:} There exists a countably infinite properly increasing chain of 
subcategories of $\SSet$
\[
  \fib_0\subsetneq\fib_1\subsetneq\fib_3\subsetneq\dots
  \subsetneq\fib_n\subsetneq\dots
\]
and corresponding countable properly decreasing chain of subcategories 
\[
  \cof_0\supsetneq\cof_1\supsetneq\cof_3\supsetneq\dots
  \supsetneq\cof_n\supsetneq\dots
\]
such that for each $n$, $\fib_n$ together with $\cof_n$ and the 
usual (topological) notion of weak equivalence provide a Quillen model 
structure on $\SSet$.  Here $\cof_0$ is the class of monomorphisms and 
$\fib_0$ that of Kan fibrations.

This phenomenon of ``variable (co)fibrations'' is quite prevalent in Quillen model categories.  Recall that Quillen in \cite{quil67} already proves the existence of two different notions of cofibration on the category of 
chain complexes of modules (with one and the same definition of weak equivalence, namely quasi-isomorphisms); on the category of simplicial diagrams, with objectwise weak equivalences, one has the cofibrations of
Bousfield--Kan~\cite{bouskan} and Heller~\cite{heller}; on cosimplicial spaces, i.e.\ cosimplicial diagrams of simplicial sets, there is yet another one due to Reedy; on the category of symmetric spectra, again at least three.

By taking products of these model categories, one sees that there is no upper bound on the cardinality of possible cofibration classes for a model category (with fixed weak equivalences, hence the same homotopy category), nor do these classes have to be ordered linearly by inclusion.  (Said differently: the identity functor may fail to be a Quillen equivalence between the same category of models, same class of weak equivalences, but with two different choices of cofibrations!)
% The different notions, alas, are not completely interchangeable, since an adjoint pair of functors may fail to be a Quillen pair (i.e.\ induce a morphism of model categories)\ with respect to one choice of cofibrations, yet do so with respect to another!
The note \cite{cofib} proves a kind of obverse to Thm.~\ref{main},
namely that in any model category satisfying mild set-theoretic
assumptions (ones satisfied by all the examples appearing
above), if one fixes the weak equivalences, and restricts cofibrations to
those that can be generated by some set, then the collection of possible
cofibration classes, partially ordered by inclusion, has least upper
bounds for any \emph{set} of elements.  It follows that the possible set-generated cofibration classes on these (combinatorial, in the sense of Jeff Smith)\ model categories all yield Quillen-equivalent model structures, the equivalence of any two arising, possibly, via a ``zig-zag'' --- not by a direct Quillen adjunction but by comparison with a third model structure.

In this paper, however, we are concerned with Thm.~\ref{main}, showing that $\SSet$ already supports infinitely many notions of cofibrations.  Topologically, they are all equivalent; the variability is due to the combinatorics of simplices.  A map belongs to the $n^\th$ exotic sense of fibration in Thm.~\ref{main},
quite simply, if it becomes a Kan fibration after $n$ iterations of Kan's
simplicial extension functor $\Ex$.  Proving that the $\fib_n$, thus defined, form part of a Quillen model
structure on $\SSet$ is straightforward.  What is surprisingly involved is proving the strict monotonicity of the inclusion $\fin_n\subsetneq\fib_{n+1}$.  It is enough to show that there is a simplicial set that becomes a Kan complex after exactly $n+1$ iterations of $\Ex$.  Using the small object argument, one can generate a fairly explicit family of simplicial sets $X$ such that $\Ex^{n+1}(X)$ is fibrant (in the ordinary sense).  The hard part is finding an $X$ among them such that $\Ex^n(X)$ is not yet fibrant.  We show, by an ad hoc path-length argument in the category of graphs, that the fibrantization (in the $n+1^{st}$ sense)\ of the $n^\th$ subdivision of the horn $\Lambda^0_2$ is such an $X$.  Many aspects of the combinatorics of iterated $\Ex$ remain delightfully mysterious; some surprising connections will be pointed out in the closing section of this paper.

As far as homotopy model theory is concerned, the role of $\Ex$ is merely one of
convenience; any adjoint of a subdivision functor satisfying basic
compatibility properties with simplices would do.  In particular, one expects that there 
are other countable decreasing chains of `axiomatic cofibrations', and the whole collection of these is no longer linearly ordered by inclusion.  Some embarrassingly natural questions remain unanswered (note that weak equivalences have been fixed throughout to be the usual ones!):

$\bullet$ Are there uncountably many distinct cofibration classes in 
$\SSet$?  Perhaps even a proper class of them? 

$\bullet$ Is there any cofibration class in $\SSet$ that is not a subclass of the 
monomorphisms?  (Equivalently, does every axiomatic class of fibrations 
include the Kan fibrations?)

From \cite{cofib} it follows that if there is merely a \emph{set} of possible cofibration classes, then there is a unique maximal one among them; but it would still be unclear whether the maximal class coincides with the monomorphisms in the case of $\SSet$.

The somewhat exotic case of non-standard cofibrations of simplicial sets has a better-known analogue in the case of sheaves of simplicial sets.  For any Grothendieck site $(\C,J)$, there exists \emph{potentially} a proper class of types of cofibrations (all of them being subclasses of monomorphisms)\ that yield a Quillen model structure on simplicial sheaves resp.\ simplicial presheaves on $(\C,J)$, weak equivalences being the usual stalkwise ones.  (See Beke~\cite{sheafi} Example 2.17, Jardine~\cite{jard06}, Isaksen~\cite{isak05}.)  Unlike the class of all monomorphisms, these intermediate cofibration classes are not functorial with respect to all geometric morphisms between toposes.  Again, I do not know whether there is actually only a set's worth of types of cofibrations, or whether the class of all monomorphisms is maximal.

Note that the intermediate cofibration classes on simplicial (pre)sheaves are ultimately due to the existence of Grothendieck topologies on the underlying category $\C$.  In the case of $\SSet$ (simplicial presheaves on the terminal category $\C$!), what this paper shows is that exotic cofibration classes arise from subdivisions of the cosimplicial object $\Delta\ra\SSet$.  I have little doubt that this phenomenon too can be exhibited on simplicial (pre)sheaves.

\section{Subdivided cofibrations}
The next proposition --- though stated for the case of simplicial sets --- is very simple, and would apply in the context of any combinatorial model category equipped with a Quillen self-adjunction (whose left adjoint part one could think of, formally, as `refining cofibrations' and its right adjoint, as `partial fibrantization').

\begin{prop}  \label{prelim}
Let $\SSet\three{\Sd^n}{\leri}{\Ex^n}\SSet$ be the $n$-fold iteration of the simplicial subdivision -- extension adjunction, and let $c_{i,k}:\Lambda^i_k\inc\Delta_k$, $k\in\N^+$, $0\leq i\leq k$ be the set of generating cofibrations for $\SSet$.  (By convention, set $\Ex^0$, $\Sd^0$ to be the identity.)  In $\SSet$, define

$\bullet$ $\cof_n$ to be the closure under pushouts, transfinite compositions and retracts of the set of morphisms $\Sd^n(c_{i,k})$

$\bullet$ $\tW$ to be the class of topological weak equivalences

$\bullet$ $\fib_n$ to be the class of morphisms $f$ such that $\Ex^n(f)$ is a Kan fibration.

\noindent
Then $\cof_n$, $\tW$ and $\fib_n$ form a Quillen model structure on $\SSet$.
\end{prop}

\begin{proof}
Consider the adjunction $\SSet\three{\Sd^n}{\leri}{\Ex^n}\SSet$ and define (for the moment)\ $\tW^{-1}$ to be the class of maps $f$ such that $\Ex^n(f)\in\tW$.  Since $\Ex^n$ preserves Kan fibrations, topological weak equivalences and arbitrary filtered colimits, it follows from the small object argument that $\cof_n$, $\tW^{-1}$ and $\fib_n$ define a Quillen model structure on $\SSet$.  (See Hirschhorn~\cite{hirsch} or Hovey~\cite{hovey} for the statement of `creating model structures by right adjoints' in the context of cofibrantly generated model categories.)

But $\tW^{-1}=\tW$ since $\Ex(f)$ is a topological weak equivalence if and only if $f$ is one; this follows from the existence of a natural inclusion $X\llra{\eta_X}\Ex(X)$ that is a weak equivalence for all $X$ (see Kan~\cite{kan57})
\[ \xymatrix{ X \ar[d]_f\ar@{>->}[r]^(.4){\eta_X} & \Ex(X) \ar[d]^{\Ex(f)} \\
              Y \ar@{>->}[r]^(.4){\eta_Y} & \Ex(Y) }\]
and the 2-of-3 property.
\end{proof}

It is now that the work begins.

\begin{thm}   \label{main}
For the model structures defined in Prop.~\ref{prelim}, one has strictly monotone chains of inclusions
\[
  \fib_0\subsetneq\fib_1\subsetneq\fib_3\subsetneq\dots
  \subsetneq\fib_n\subsetneq\dots
\]
resp.\
\[
  \cof_0\supsetneq\cof_1\supsetneq\cof_3\supsetneq\dots
  \supsetneq\cof_n\supsetneq\dots
\]
\end{thm}

\emph{Outline of the proof.} Since $\Ex$ preserves Kan fibrations, the inclusion $\fib_n\subseteq\fib_{n+1}$ is automatic, and that implies $\cof_n\supseteq\cof_{n+1}$.  The strictness follows from

\begin{prop}   \label{decr}
For any $n\in\N$, there exists a simplicial set $X$ such that $\Ex^n(X)$
does not satisfy the Kan extension condition, but $\Ex^{n+1}(X)$
does.
\end{prop}

The proof is preceded by two lemmas.  The first one states, roughly, that in the $n^\th$ barycentric subdivision of a triangle, pairs of points on the boundary cannot be connected by interior paths shorter than $2^n$.  (This will be responsible for non-injectivity of a certain graph with respect to certain graph maps.)  The second lemma states an analogue of this for the $n^\th$ simplicial subdivision of the simplex $\Delta_n$.  We then exhibit the required counterexample $X$: it is $R_\infty(\Sd^n\Lambda^0_2)$, where $R_\infty$ is the canonical fibrantization functor for the model structure $\fib_n$, and $\Lambda^0_2$ is $\Delta_2$ minus its (non-degenerate)\ 2-simplex and $0^\th$ face.  (See Conj.~\ref{nervie} for another guess at where counterexamples may come from.)

The first lemma, for the sake of visual simplicity, will be stated for `barycentric subdivisions' in the classical sense of simplicial complexes.

\begin{lemma}  \label{2d}
Let $x$ and $y$ be vertices of the $n^\th$ barycentric subdivision of a triangle with vertices $A$, $B$, $C$.  Suppose $x$ lies on the side $AB$ and $y$ on the side $AC$ of the triangle, $x\neq A$ and $y\neq A$.  Let $p$ be an edge path connecting $A$ and $B$.  Suppose $p$ does not pass through the vertex $A$.  Then $p$ contains at least $2^n$ edges.
\end{lemma}

This is an example of the statement for $n=2$:

\[\begin{xy}
<1cm,0cm>:
(8,0)**@{.}
;(4,6)**@{.}
;(0,0)**@{.}
;(6,3)**@{.}
,(8,0);(2,3)**@{.}
,(4,0);(4,6)**@{.}
,(0,0);(6,3)**@{.}
;(4,4)**@{.};(2,3)**@{.};(2,1)**@{-}
;(2.666,.666)**@{-};(4,0)**@{.};(6,1)**@{.};(6,3)**@{.}
,(0,0);(4,1)**@{.};(8,0)**@{.};(5,2.5)**@{.}
;(4,6)**@{.};(3,2.5)**@{.};(0,0)**@{.}
,(2,0);(2.666,.666)**@{-};(4,2)**@{.};(6,0)**@{.}
,(7,1.5);(4,2)**@{.};(5,4.5)**@{.}
,(3,4.5);(4,2)**@{.};(1,1.5)**@{.}
,(8.3,-.2)*{B}
,(-.3,-.2)*{A}
,(4,6.2)*{C}
,(2,-.2)*{x}
,(1.7,3)*{y}
\end{xy}\]   

\begin{center}
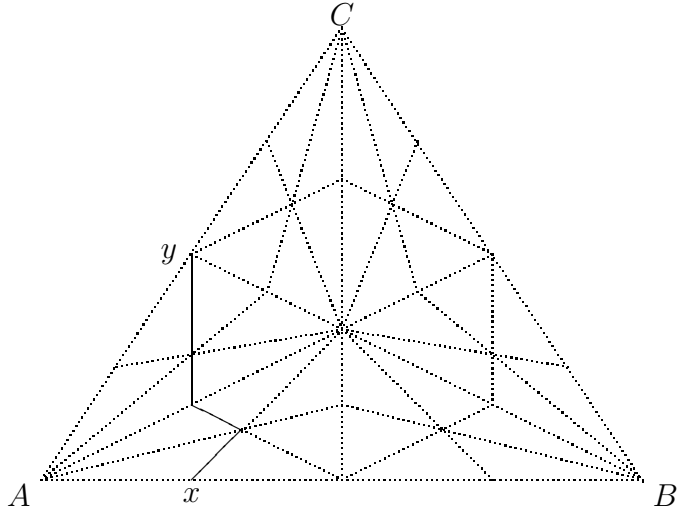
{\small \textsc{Fig.~1.} To get from side $AB$ to side $AC$, avoiding the 
vertex $A$,\\you need an edge path of length at least $2^2$ in the twice-subdivided triangle $ABC$.}
\end{center}

The proof of this lemma (which is an inductive partitioning argument)\ is postponed.  In what follows, we return to the world of simplicial sets (with degeneracies)\ and simplicial subdivisions.  For vertices $x$, $y$ of a simplicial set, write $d(x,y)$ for their edge distance, that is to say, the least length of a \emph{possibly ``zig-zag'' edge path} connecting them.  (All simplicial sets considered below will be connected.)

If $x,y$ are vertices of $X$, and $X\llra{f}Y$ is a map of simplicial sets, note that
\[
        d(x,y)\geq d(f(x),f(y)).
\]
\begin{lemma}  \label{3d}
Write $\partial\Delta_k$ for the boundary of the standard $k$-simplex, $k>1$.
Let $x$ and $y$ be vertices of $\Sd^n(\partial\Delta_k)$, thought
of as simplicial subset of $\Sd^n(\Delta_k)$.  If $d(x,y) < 2^n$ in
$\Sd^n(\Delta_k)$, then the distance of $x$ and $y$ in
$\Sd^n(\partial\Delta_k)$ equals their distance in
$\Sd^n(\Delta_k)$.
\end{lemma}

\begin{proof}
(a) Suppose there is a top-dimensional face $i:\Delta_{k-1}\inc\Delta_k$ of our $k$-simplex such that $\Sd^n(\Delta_{k-1})$ contains both $x$ and $y$.  There is a retraction $r:\Delta_k\ra\Delta_{k-1}$ in $\SSet$ (a degeneracy `dual' to $i$), whence a retraction $\Sd^n(r)$; by the above remark, the distance
of $x$ and $y$ in $\Sd^n(\partial\Delta_k)$ then cannot be greater than their distance in $\Sd^n(\Delta_k)$.

(b) If no face of $\Delta_k$ contains both $x$ and $y$, then, without loss of generality, assume
that $x$ lies on the face opposite the vertex $[0]$, $y$ lies on the face
opposite the vertex $[1]$, and neither lies on the intersection of these
faces, the (codimension 2) face $\F$ with vertices $[2],[3],\dots,[k]$.  
Consider a distance-minimizing path $p$ in $\Sd^n(\Delta_k)$ between
$x$ and $y$.  If $p$ contains a vertex $F$ on the subdivided face
$\F$, then the argument of part (a) can be applied separately to the paths $XF$
and $FY$ to deduce that a distance-minimizing edge path between $x$ and
$y$ can proceed on $\Sd^n(\partial\Delta_k)$, as claimed.

(c) The missing case is when the distance-minimizing path $p$ avoids $\F$.  We show that any such
path must be of length $2^n$ at least, contradicting our assumption
that $d(x,y) < 2^n$.

Consider the simplicial collapsing map $\Delta_k\llra{c}\Delta_2$ corresponding to the monotone map that sends $[0]$ to $[0]$, $[1]$ to $[1]$, and $[i]$ to $[2]$ for $i\geq 2$.  Under the map $\Sd^n(p))$, $x$ and $y$ are sent into vertices of $\Sd^n(\Delta_2)$, $x$ lying on
the side opposite the vertex $[0]$, $y$ lying on the side opposite the
vertex $[1]$, and $p$ will become an edge path connecting them that avoids
the vertex $[2]$.  By Lemma~\ref{2d}, $\Sd^n(p))$ has at least length $2^n$ (note that the presence of degenerate edges does not change path distances), whence so does $p$.
\end{proof}

For any simplicial set $U$, define $R_\infty(U)$ to be the colimit of the chain
\[
   U=:R_0(U) \ra R_1(U) \ra R_2(U) \ra R_3(U) \ra \dots
\]
where $R_{j+1}(U)$ arises from $R_j(U)$ by pushing on all $n+1$-times subdivided horn filling
conditions
\[ \xymatrix{ \Sd^{n+1}(\Lambda^i_k) \ar[r]\ar@{>->}[d] & R_j(U) \\
              \Sd^{n+1}(\Delta_k)   }   \]
that exist at that stage.  By Quillen's small object argument, $R_\infty(U)$ has the right lifting property with respect to the set of maps $\Sd^{n+1}(\Lambda^i_k)\ra\Sd^{n+1}(\Delta_k)$.  Adjointly, $\Ex^{n+1}(R_\infty(U))$ is a Kan complex.  Set $U=\Sd^n(\Lambda^0_2)$.  We will exhibit a specific lifting problem with respect to an $n$-times subdivided horn inclusion that $X=R_\infty(\Sd^n(\Lambda^0_2))$ \emph{fails}; that is to say, $\Ex^n(X)$ is not a Kan complex.

The lifting problem will be
\[ \xymatrix{
     \Sd^n(\Lambda^0_2) \ar[rr]^(.45){\text{canonical}}\ar@{>->}[d]^i &&
                                     R_\infty(\Sd^n(\Lambda^0_2)) \\
            \Sd^n(\Delta_2) \ar@{.>}[urr]_{(?)}   } \]
$\Sd^n(\Lambda^0_2)$ is precisely a zig-zag of length $2^{n+1}$.  
Call its extreme vertices $x$ and $y$.  ($x$ and $y$ are thus the vertices 
of $\Delta_2$ that bound the edge missing in $\Lambda^0_2$.)  Note that 
$d(i(x),i(y))=2^n$ in $\Sd^n(\Delta_2)$.  If a lift $(?)$ existed, 
then it would have to exist into $R_j(\Sd^n(\Lambda^0_2))$ for 
some finite $j$ already, since $\Sd^n(\Delta_2)$ is (simplicially)\ finite.  So, 
letting $r_j$ denote the canonical map
$\Sd^n(\Lambda^0_2)\ra R_j(\Sd^n(\Lambda^0_2))$, to prove the 
impossibility of a lift, it suffices to prove

\qquad $d(r_j(x),r_j(y)) = 2^{n+1}$ in $R_j(\Sd^n(\Lambda^0_2))$ for all $j\geq 0$.

\noindent
This is true for $j=0$; now use induction.  $R_{j+1}(\Sd^n(\Lambda^0_2))$ arises from 
$R_j(\Sd^n(\Lambda^0_2))$ via simultaneous pushouts of the type
\[ \xymatrix{ \Sd^{n+1}(\Lambda^i_k) \ar[r]\ar@{>->}[d]
                      & R_j(\Sd^n(\Lambda^0_2)) \ar[d] \\
 \Sd^{n+1}(\Delta_k) \ar[r] & R_{j+1}(\Sd^n(\Lambda^0_2)) }\]

Let $a,b$ be any vertices of $R_j(\Sd^n(\Lambda^0_2))$.  In $R_{j+1}(\Sd^n(\Lambda^0_2))$, possibly new paths have been pushed on that connect $a$ and $b$, but by Lemma~\ref{3d}, if $d<2^{n+1}$, paths of length $d$ are attached only between $a,b$ whose distance in $R_j(\Sd^n(\Lambda^0_2))$ is $d$ or less.  Therefore, if $d(a,b)\leq 2^{n+1}$, the distance of $a$ and $b$ in $R_j(\Sd^n(\Lambda^0_2))$ equals their distance in
$R_{j+1}(\Sd^n(\Lambda^0_2))$.  In particular, by the induction hypothesis, $d(r_j(x),r_j(y)) = 2^{n+1} = d(r_{j+1}(x),r_{j+1}(y))$.

To finish the proof of Prop.~\ref{decr}, we still need to prove Lemma~\ref{2d}.  Let us return to the language of planar figures.  By induction on $n$, we will show that the $6^n$ triangles in the $n^\th$ barycentric subdivision of $ABC$ can be assigned into $2^n$ disjoint classes (which we will call `rays' and label with the integers from $1$ through $2^n$)\ such that
\begin{itemize}
\item[(1)] Side $AB$ (other than the vertex $A$ itself)\ lies on ray 1; side $AC$ (other than the vertex $A$)\ lies on ray $2^n$.
\item[(2)] Let $T$ be one of the $6^n$ small triangles.  Suppose $T$ belongs to ray $i$ and does not contain the vertex $A$.  Then either $(2a)$ one edge of $T$ lies on the common boundary of ray $i$ and ray $i+1$ (for some $1\leq i\leq 2^n$)\ and its opposite vertex lies on the common boundary of ray $i$ and ray $i-1$ or $(2b)$ one edge of $T$ lies on the common boundary of ray $i$ and ray $i-1$ (for some $1\leq i\leq 2^n$)\ and its opposite vertex lies on the common boundary of ray $i$ and ray $i+1$.  The interior of the other two edges of $T$, in both cases, will belong to ray $i$.
\end{itemize}
(To avoid having to state separate cases for $i=0$ and $i=2^n$, let us agree that the side $AB$ belongs to ray 0, and side $AC$ belongs to ray $2^n+1$.)

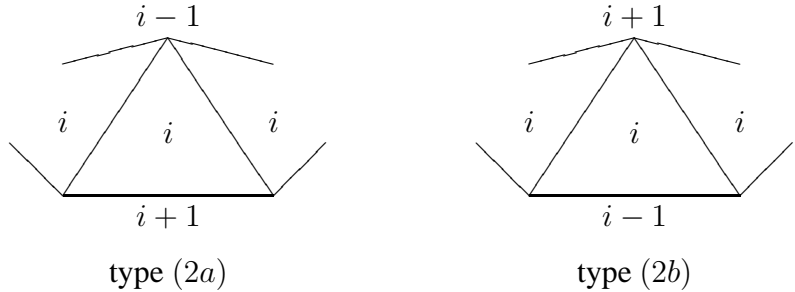
\[\begin{xy}
<7mm,0mm>:
(4,0)**@{-}
;(2,3)**@{-}
;(0,0)**@{-}
,(2,-.4)*{i+1}
,(2,-1.5)*{\textup{type\ }(2a)}
,(2,3.4)*{i-1}
,(2,1.2)*{i}
,(0,1.4)*{i}
,(4,1.4)*{i}
,(4,0);(5,1)**@{-}
,(2,3);(4,2.5)**@{-}
,(2,3);(0,2.5)**@{-}
,(0,0);(-1,1)**@{-}
\end{xy}
\hspace{2cm}
\begin{xy}
<7mm,0mm>:
(4,0)**@{-}
;(2,3)**@{-}
;(0,0)**@{-}
,(2,-.4)*{i-1}
,(2,-1.5)*{\textup{type\ }(2b)}
,(2,3.4)*{i+1}
,(2,1.2)*{i}
,(0,1.4)*{i}
,(4,1.4)*{i}
,(4,0);(5,1)**@{-}
,(2,3);(4,2.5)**@{-}
,(2,3);(0,2.5)**@{-}
,(0,0);(-1,1)**@{-}
\end{xy}\]

From (1)\ and (2)\ it follows that an interior edge of the subdivided triangle, if it does not contain the vertex $A$, either lies on the common boundary of ray $i$ and ray $i+1$ (for some $1\leq i\leq 2^n$), or spans ray $i$ (so that one of its endpoints belongs to ray $i-1$ and ray $i$, and the other endpoint belongs to ray $i$ and ray $i+1$).  To get from point $x$ on the side $AB$ to point $y$ on the side $AC$, avoiding vertex $A$, a path must cross all $2^n$ rays, so must contain at least $2^n$ edges, as claimed.

Here are the rays for $n=1$ and $n=2$.

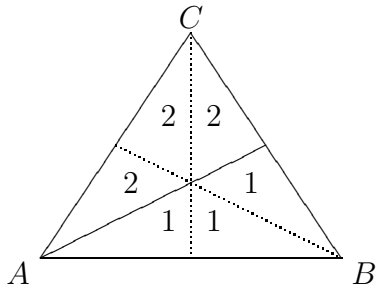
\[\begin{xy}
<1cm,0cm>:
(4,0)**@{-}
;(2,3)**@{-}
;(0,0)**@{-}
;(3,1.5)**@{-}
,(4,0);(1,1.5)**@{.}
,(2,0);(2,3)**@{.}
,(4.3,-.2)*{B}
,(-.3,-.2)*{A}
,(2,3.2)*{C}
,(1.7,.5)*{1}
,(2.3,.5)*{1}
,(2.8,1)*{1}
,(1.2,1)*{2}
,(1.7,1.9)*{2}
,(2.3,1.9)*{2}
\end{xy}\]

\begin{center}
{\small \textsc{Fig.~2.} The partitioning of $\Sd^1\Delta_2$.}
\end{center}

\[\begin{xy}
<1cm,0cm>:
(8,0)**@{-}
;(4,6)**@{-}
;(0,0)**@{-}
;(6,3)**@{-}
,(8,0);(2,3)**@{.}
,(4,0);(4,6)**@{.}
,(0,0);(4,1)**@{-};(5.333,.666)**@{-}
;(6,1)**@{-};(6,1.666)**@{-};(7,1.5)**@{-}
,(2,0);(4,2)**@{.};(1,1.5)**@{.}
,(4,0);(2,1)**@{.}
,(4,0);(5.333,.666)**@{.};(8,0)**@{.}
,(6,0);(4,2)**@{.};(6,1.666)**@{.};(6,3)**@{.}
;(4.666,3.666)**@{.};(4,2)**@{.};(3,4.5)**@{.}
,(8,0);(5,2.5)**@{.};(4,6)**@{.};(3.333,3.666)**@{.}
;(2,3)**@{.};(2,1)**@{.}
,(5,4.5);(4.666,3.666)**@{-}
;(4,4)**@{-};(3.333,3.666)**@{-};(3,2.5)**@{-};(0,0)**@{-}
,(3.6,.5)*{1}
,(5.5,1.5)*{2}
,(4.2,3)*{3}
,(2.9,3.9)*{4}
,(8.3,-.2)*{B}
,(-.3,-.2)*{A}
,(4,6.2)*{C}
\end{xy}\]

\begin{center}
{\small \textsc{Fig.~3.} The partitioning of $\Sd^2\Delta_2$.\\
Only one triangle in each contiguous region is marked with its number $i$.}
\end{center}

In general, the partitions are defined by induction.  Let $T$ be a triangle of $\Sd^n\Delta_2$, not containing the vertex $A$, and of the type that was denoted $(2a)$ above.  Its subdivisions will then be assigned numbers

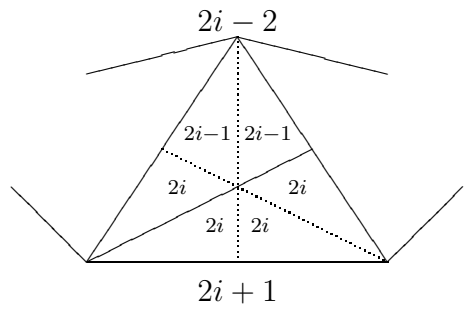
\[\begin{xy}
<1cm,0cm>:
(4,0)**@{-}
;(2,3)**@{-}
;(0,0)**@{-}
;(3,1.5)**@{-}
,(4,0);(1,1.5)**@{.}
,(2,0);(2,3)**@{.}
,(2,3.2)*{2i-2}
,(1.7,.5)*{\msm{2i}}
,(2.3,.5)*{\msm{2i}}
,(2.8,1)*{\msm{2i}}  
,(1.2,1)*{\msm{2i}}
,(1.6,1.7)*{\msm{2i-1}}
,(2.4,1.7)*{\msm{2i-1}}
,(2,-.4)*{2i+1}
,(4,0);(5,1)**@{-}
,(2,3);(4,2.5)**@{-}
,(2,3);(0,2.5)**@{-}
,(0,0);(-1,1)**@{-}
\end{xy}\]

If $T$ is of type $(2b)$, its subdivisions will be labeled

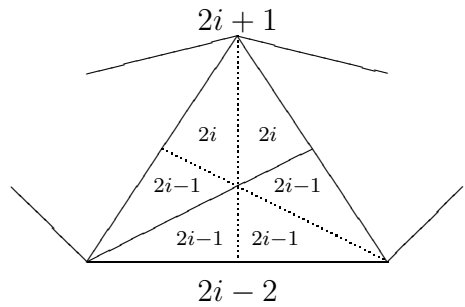
\[\begin{xy}
<1cm,0cm>:
(4,0)**@{-}
;(2,3)**@{-}
;(0,0)**@{-}
;(3,1.5)**@{-}
,(4,0);(1,1.5)**@{.}
,(2,0);(2,3)**@{.}
,(2,3.2)*{2i+1}
,(1.5,.3)*{\msm{2i-1}}
,(2.5,.3)*{\msm{2i-1}}
,(2.8,1)*{\msm{2i-1}}  
,(1.2,1)*{\msm{2i-1}}
,(1.6,1.7)*{\msm{2i}}
,(2.4,1.7)*{\msm{2i}}
,(2,-.4)*{2i-2}
,(4,0);(5,1)**@{-}
,(2,3);(4,2.5)**@{-}
,(2,3);(0,2.5)**@{-}
,(0,0);(-1,1)**@{-}
\end{xy}\]

\noindent
and the induction hypotheses are satisfied.  (It is worthwhile to iterate
the construction and observe the `fractal boundaries' of the rays, and the self-similarity of the local patterns arising.)

As for the triangles in $\Sd^n\Delta_2$ that contain the vertex $A$, forming a fan around $A$, they are numbered consecutively from 1 (at side $AB$)\ to $2^n$ (at side $AC$); this is compatible with the subdivisions of type $(2a)$ and $(2b)$.

This, then, finishes the proof of Lemma~\ref{2d}, and also of Prop.~\ref{decr}, so of the main theorem. \qed

\noindent
One can show (see Prop.~\ref{latch} below)\ that the standard simplices belong to $\fib_n$ for $n>0$.  ($\Delta_k$ itself is a Kan complex only for $k=0$.)  On the other hand, for $n>0$ it will no longer be true that every simplicial set is cofibrant.

\section{On the way to fibrancy}
The path length counterexample was quite artificial, and it may be of combinatorial interest to understand how other families of simplicial sets --- for example, nerves $N\C$ of categories $\C$ --- become fibrant.  The following two facts are recalled as `teasers'.

\begin{prop}
For a small category $\C$, $N\C$ is a Kan complex if and only if $\C$ is a groupoid.
\end{prop}

This is classical; a proof can be found in e.g.\ Lee~\cite{lee72}.

\begin{prop}   \label{latch}
For a small category $\C$, $\Ex(N\C)$ is a Kan complex if and only if $\C$ (thought of as a full subcategory of itself)\ possesses a left calculus of fractions.
\end{prop}

This is Latch--Thomason--Wilson~\cite{ltw}, remark 5.8.  Note that a category being a groupoid amounts to injectivity with respect to two functors in the category of (small) categories: these ensure the possibility of left and right ``division''.  Similarly, the property of a category ``possessing a left calculus of fractions with respect to itself'' amounts to injectivity with respect to the following two functors between finite diagrams:
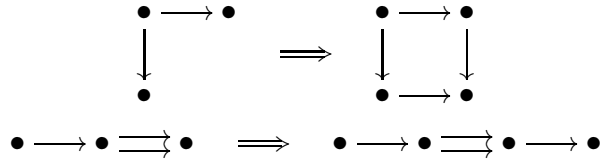
\[ \xymatrixcolsep{.5pc}\xymatrixrowsep{.5pc} \xymatrix{
\bullet \ar[rr]\ar[dd] && \bullet &&&& \bullet \ar[rr]\ar[dd] && \bullet
\ar[dd] \\ &&& \ar@{=>}[rr] &&  \\
\bullet &&&&&& \bullet \ar[rr] && \bullet } \]
\[ \xymatrixcolsep{.5pc}\xymatrixrowsep{.5pc} \xymatrix{
\bullet \ar[rr] &&\bullet \ar@<.5ex>[rr]\ar@<-.5ex>[rr] &&\bullet
& \ar@{=>}[rr] &&{} & \bullet \ar[rr] &&
\bullet \ar@<.5ex>[rr]\ar@<-.5ex>[rr] &&\bullet \ar[rr] &&\bullet } \]
(Composition rules for arrows are omitted, but see Gabriel--Zisman~\cite{gabzis}.)

One can think of $\C$ ``possessing a left calculus of fractions with respect to all morphisms'' as an approximation to $\C$ being a groupoid; the morphisms of $\C[\C^{-1}]$, while not actual arrows, are representable by equivalence classes of zig-zags of length 2.  If one takes $\C$ to have a single object (i.e.\ to be a monoid), then for it to have a left calculus of fractions means that it satisfies the left Ore conditions; in a certain way, it is close to being a group.  It is tempting to make the following

\begin{conj}
For each $n$, there exists a \emph{finite} collection $I$ of functors between finitely presentable categories such that $\Ex^n(N\C)$ is a Kan complex if and only if $C$ is injective with respect to $I$.
\end{conj}

\begin{conj}    \label{nervie}
For any $n$, there exist finitely presentable categories $\C$ (even monoids)\ such that $\Ex^n(N\C)$ is not a Kan complex, but $\Ex^{n+1}(N\C)$ is.
\end{conj}

$\Ex^n(N\C)$ is a Kan complex if and only if the category $\C$ is injective with respect to a certain family of maps between finite posets --- namely, $\cat\big(\Sd^n(\Lambda_k^i\inc\Delta_k)\big)$, where $\cat$ is the left adjoint of the nerve functor --- but it is far from obvious that this is also equivalent to a \emph{finite} collection of injectivity conditions.  (One can show that the notion of Kan fibration cannot be axiomatized by finitely many injectivity conditions between finite simplicial sets; in fact, not even by finitely many first-order axioms in the language of simplicial sets.)

Let us now consider the homotopy types of categories
whose nerves become Kan complexes after finitely many iterations of $\Ex$.  If $\C$ is
a groupoid, then $N\C$ is homotopy equivalent to the disjoint union of
Eilenberg--MacLane spaces $K(\pi,1)$ corresponding to its vertex groups.  
By a theorem of Dwyer and Kan~\cite{dk80b}, if $\C$ possesses a left or
right (more generally, homotopy left or right) calculus of fractions with
respect to all its morphisms, then the localization map $\C\ra\C[\C^{-1}]$
induces a weak equivalence on nerves.  By putting this and
Prop.~\ref{latch} together, one sees that the $n=0$ and $n=1$ cases of the following conjecture hold:

\begin{conj}
Suppose that for a small category $\C$, $\Ex^n(N\C)$ is a Kan complex for 
a finite $n$.  Then $N\C$ is weakly equivalent to a disjoint union of Eilenberg--MacLane spaces (i.e.\ has vanishing homotopy groups above dimension 1).
\end{conj}

This conjecture is, of course, rather daring.  The underlying intuition is that as $n$ increases, one has a progressive weakening of the algebraic notion of group(oid) --- the $n^\th$ such weakening being that $\C$ is such that $\Ex^n(N\C)$ is a Kan complex --- but all these notions are special cases of ``groupoids up to homotopy''.  (This is rather in line with the philosophy of certain approaches to higher categories.)

Let us leave these conjectures now.  At this stage, the reader has no doubt already called to mind the work of
Thomason~\cite{thom80}; he proves that the categorification-nerve
adjunction $\Cat\three{\cat}{\leri}{N}\SSet$ does not create a Quillen
model structure on $\Cat$ from the one on $\SSet$ we denoted
$\lc\cof_n,\tW,\fib_n\rc$ for $n=0$ and $n=1$, but does create
one for $n=2$.  It follows that it creates one from
$\lc\cof_n,\tW,\fib_n\rc$ for any $n\geq 2$; but it does not follow that the fibrancy classes of these model structures on $\Cat$ are distinct.  Conj.~\ref{nervie} implies that they are.  (Perhaps this can also be proved by the methods of this paper, looking carefully at composability of edges in $\Sd^n\Delta_k$.)  This situation is typical when one transports the ``subdivided cofibrations'' model structures on $\SSet$ (as arising either from Kan's $\Sd$ or from another simplicial subdivision functor)\ to groupoids, small categories, simplicial universal algebras, etc., or sheafifies them~\cite{sheafi}: one needs to check whether new fibrations actually arise.  An example when this does \emph{not} happen is across the adjunction
\[
    \Top\three{|-|}{\leri}{\Sing}\SSet
\]
owing to the fact that the geometric realization of a subdivided simplex 
is homeomorphic to the original.  Maybe (compactly generated)\ topological 
spaces and weak equivalences possess an extremal fibration class.

\bibliographystyle{plain}

\end{document}